\documentstyle[12pt]{article}

\begin{document}

\date{}

\title{\Large\bf Noncommutativity and Discrete Physics}

\author{Louis H. Kauffman \\
  Department of Mathematics, Statistics and Computer Science \\
  University of Illinois at Chicago \\
  851 South Morgan Street,
  Chicago, IL, 60607-7045}

 \maketitle

 \thispagestyle{empty}

 \section{Introduction}

The purpose of this paper is to present an introduction to a point of view
for discrete foundations of physics. In taking a discrete stance, we find
that the initial expression of physical theory must occur in a context of
noncommutative algebra and noncommutative vector analysis. In this way the
formalism of quantum mechanics occurs first, but not necessarily with the
usual interpretations. By following this line carefully we can show how the
outlines of the well-known commutative forms of physical theory arise first
in noncommutative form. This much, the present paper will make clear with
specific examples and mathematical formulations. The exact relation of
commutative and noncommutative theories raises a host of problems.
\vspace{3mm}

In the first section of this paper we discuss  the properties of the
noncommutative discrete calculus that underlies our work. The section ends
with the consequences in our framework for a particle whose position and
momentum commutator is equated to a (noncommutative) metric field. In the
next section we discuss how our discrete stance leads to an inversion of
the usual Dirac maxim "replace Poisson brackets with commutators". If we
replace commutators with Poisson brackets that obey a Leibniz rule
satisfied by our commutators, then the dynamical variables will obey
Hamilton's equations. Thus we can take Hamilton's equations as the natural
classicization of our theory. The next section shows how the
noncommutativity is neccessary in this approach and shows how certain
representations of the theory lead to chaotic dynamics.  The next section
discusses the relationship of the discrete ordered calculus with
$q$-deformations and quantum groups. We show that in a quantum group with a
special grouplike element representing the square of the antipode, there is
a representation of the discrete ordered calculus. In this calculus on a
quantum group the square of the antipode represents one tick of the clock.
Then follows a section on networks and discrete spacetime. This section is
a general exposition of ideas related to spin networks and topological
quantum field theory. It is our speculation that the approaches to discrete
physics inherent in discrete calculus and in topological field theory are
deeply interrelated.  At the end of this section we outline this
relationship in the case of a recent model for quantum gravity due to Louis
Crane.
\vspace{3mm}

\noindent
{\bf Acknowledgement.}    It gives the author pleasure to thank the
National Science Foundation for support of this research under
NSF Grant  DMS -2528707.
 \vspace{3mm}

\section{Discrete Ordered Calculus}

Consider successive measurements of position and velocity. In measurement
of position, no time step is required. In measuring velocity, we need
positions of two neighboring instants of time.
\vspace{3mm}

Thinking discretely, let us assume that the particle has positions
$$X,X',X'',....$$
at successive moments of time. Discrete unit time steps are indicated by
the primes appended to the $X$. A general point in the time series at time
$t$ will be denoted by $X^{t}$.  By convention let the time step between
successive points in the series be equal to 1 :  $$\Delta t = 1$$.
Then  we can define the velocity at time $t$ by the formula:
$$v(t) = X^{t+1} - X^{t}.$$
More generally, if $X$ denotes the position at a given time, then $X'-X$
denotes the velocity {\em at that time}, where the phrase "at that time"
must involve the next time as well. In a discrete context there is no
notion of instantaneous velocity.
\vspace{3mm}

Measure position, and you find $X$. Then measure velocity, and you get $X'
- X$. Now measure position, and you get $X'$ because the timehas shifted
to the next time in order to allow the velocity measurement. In order to
measure velocity the position is necessarily shifted to its value at the
next time step. In this sense, position and velocity measurements cannot
commute in a discrete framework.
\vspace{3mm}

Our project is to take this basic noncommutativity at face value and follow
out its consequences. To this end we will formulate a calculus of finite
differences that takes the order of observations into account. This
formalization is explained below. Remarkably, the resulting calculus is
actually a discrete version of time evolution in standard quantum
mechanics.
\vspace{3mm}

We begin by recalling the usual derivative in the calculus of finite
differences, generalised to a (possibly) noncommutative context.
\vspace{3mm}

\noindent
{\bf Definition.} Let $$dX = X'-X$$ define the finite difference derivative
of a variable $X$ whose successive values in discrete time are
$$X,X',X'',....$$  This $dX$ is a classical derivative in the calculus of
finite differences. It is still defined even if the quantities elements of
the time series are in a noncommutative algebra.  We shall assume that the
values of the time series are in a possibly noncommutative ring $R$ with
unit. (Thus the values could be real numbers, complex numbers, matrices,
linear operators on a Hilbert space, or elements of an appropriate abstract
algebra.) This means that for every element $A$ of the ring $R$ there is a
well-defined successor element $A'$, the next term in the time series. It
is convenient to assume that the ring itself has this temporal structure.
In practice, one is concerned with a particular time series and not the
structure of the entire ring.  Moreover, we shall assume that the next-time
operator distributes over both addition and multiplication in the sense
that
$$(A+B)' = A' + B'$$  and $$(AB)' = A'B'.$$
An element $c$  of the ring $R$ is said to be a $constant$ if $c' = c.$
\vspace{3mm}

\noindent
{\bf Lemma.} $$d(XY)= X'd(Y) + d(X)Y.$$
\vspace{3mm}

\noindent
{\bf Proof.} $$d(XY) = X'Y'-XY$$
$$=X'Y'-X'Y+X'Y-XY$$
$$=X'(Y'-Y) +(X'-X)Y$$
$$=X'd(Y) + d(X)Y.$$
\vspace{3mm}

This formula is {\em different} from the usual formula in Newtonian
calculus by the time shift of $X$ to $X'$ in the first term. We now correct
this discrepancy in the calculus of finite differences by taking a {\em
new} derivative $D$ as an {\em instruction to shift the time to the left of
the operator $D$.} That is, we take $XD(Y)$ quite literally as an
instruction to {\em first find $dY$ and then find the value of $X.$} In
order to find $dY$ the clock must advance one notch. Therefore $X$ has
advanced to $X'$ and we have that the evaluation of $XD(Y)$ is
$$X'(Y'-Y).$$
\vspace{3mm}

In order to keep track of this noncommutative time-shifting, we will write
$$DX= J(X'-X)$$ where the element $J$ is a special time-shift operator
satisfying $$ZJ = JZ'$$ for any $Z$ in the ring $R$. The time-shifter, $J$,
acts to automatically evaluate expressions in the resulting noncommutative
calculus of finite differences. We call this calculus $DOC$ (for discrete
ordered calculus). Note that $J$ formalizes the operational ordering
inherent in our initial discussion of velocity and position measurements.
An operator containing $J$  causes a time shift in the variables or
operators to the left of $J$ in the sequence order.
\vspace{3mm}

Formally, we extend the ring of values $R$ (see the definition of $d$
above) by adding a new symbol $J$ with the property that $AJ = JA'$ for
every $A$ in $R$. It is assumed that the extended ring $R$ is associative
and satisfies the distributive law so that $J(A+B) = JA + JB$ and $J(AB) =
(JA)B$ for all $A$ and $B$ in the ring. We also assume that $J$ itself is a
constant in the sense that $J' = J$.
\vspace{3mm}

The key result in $DOC$  is the following adjusted difference formula:
\vspace{3mm}

\noindent
{\bf Lemma 2.} $$D(XY) = XD(Y) + D(Y)X.$$
\vspace{3mm}

\noindent
{\bf Proof.} $$D(XY)$$
$$ = J(X'Y'-XY)$$
$$= J(X'Y'-X'Y +X'Y-XY)$$
$$=J(X'(Y'-Y)+ (X'-X)Y$$
$$=JX'(Y'-Y) + J(X'-X)Y$$
$$=XJ(Y'-Y) + J(X'-X)Y$$
$$=XD(Y) + D(X)Y.$$
\vspace{3mm}

The upshot is that DOC behaves formally like infinitesimal calculus and can
be used as a calculus in this version of discrete physics. In \cite{KN:QEM}
Pierre Noyes and the author use this foundation to build a derivation of  a
noncommutative version of electromagnetism.  Another version of this
derivation can be found in \cite{Twist}.  In both cases the derivation is a
translation to this context of the well-known Feynman-Dyson derivation of
electromagnetic formalism from commutation relations of position and
velocity.
\vspace{3mm}

 Note that the
definition of the derivative in DOC  is actually a commutator: $$DX =
J(X'-X) = JX' - JX = XJ -JX = [X,J].$$  The operator $J$ can be regarded as
a discretised time-evolution operator in the Heisenberg formulation of
quantum mechanics. In fact we can write formally that $$X' = J^{-1}XJ$$
since $JX' = XJ$ (assuming for this interpretation that the operator $J$ is
invertible). Putting the time variable back into the equation, we get the
evolution
$$X^{t+ \Delta t} = J^{-1}X^{t}J.$$  This aspect can be compared to the
formalism of Alain Connes' theory of noncommutative geometry \cite{Connes}.
\vspace{3mm}

In Connes' theory there is a notion of quantized differential that takes
the form (in his language)  $de = [F,e]$ where
$F$ is a bounded operator on a Hilbert space $H$ and $[e]$ is a class in
the $K-$ theory of a certain algebra $A$ acting on the Hilbert space. In
this context Connes' quantized calculus is used to obtain a wide range of
connections with various aspects of physics, including a new view of the
standard model for fundamental particles. Our approach to aspects of the
formalism of the $DOC$ quantized calculus may fit into the context of
Connes' theory.  This is a topic that derserves futher investigation. In
this paper, and in our previous work we have used the most elementary
non-commutative algebraic tools to obtain our results. It is our hope that
these results will fit into more complex contexts that are directly related
to both theory and measurement.
\vspace{3mm}

In the discrete ordered calculus,  $X$ and $DX$ have no reason to commute:
$$[X,DX] = XJ(X'-X) -J(X'-X)X = J(X'(X'-X) -(X'-X)X)$$ Hence $$[X,DX] =
J(X'X'-2X'X + XX).$$
This is non-zero even in the case where X and X' commute with one another.
Consequently, we can consider physical laws in the form $$[X^{i}, DX^{j}] =
g^{ij}$$
where $g^{ij}$ is a function that is suitable to the given application. In
\cite{KN:QEM} we show how the formalism of electromagnetism arises when
$g^{ij}$
is $\delta^{ij}$, the Kronecker delta. In \cite{KN:DG} we show how
the general case corresponds to a "particle" moving in a noncommutative
gauge field coupled with geodesic motion relative to the Levi-Civita
connection associated with the $g^{ij}.$ This result can be used to place
the work of Tanimura \cite{Tanimura} in a discrete context.
\vspace{3mm}

It should be emphasized that all physics that we derive in this way is
formulated in a context of noncommutative operators and variables. We do
not derive electromagnetism, but rather a noncommutative analog.  It is not
yet clear just what these noncommutative physical theories really mean. Our
initial idealisation of measurement is not the only model for measurement
that corresponds to actual observations. Certainly the idea that we can
measure time  in a way that "steps between the steps of time" is an
idealisation. It happens to be an idealisation that fits a model of the
universe as a cellular automaton. In a cellular automaton an observation
is what an operator of the automaton might be able to do. It is not
necessarily what the "inhabitants" of the automaton can perform. Here is
the crux of the matter. The inhabitants can have only limited observations
of the running of the automaton, due to the fact that they themselves are
processes running on the automaton. I believe that the theories we build on
the basis of DOC are theories {\em about} the structure ofthese automata.
They will eventually lead to theories of what the processes that run on
such automata can observe. It is quite possible that the well known
phenomena of quantum mechanics will arise naturally in such a context.
These points of view should be compared with \cite{Fkin}.
\vspace{3mm}

In order to illustrate these methods, I will show part of the calculations
related to $$[X^{i}, \dot{X^{j}}] = g^{ij}.$$ Here $\dot{X^{j}}$ is
shorthand for$DX^{j}.$ Along with this commutator equation, we will assume
that $$[X^{i}, X^{j}] = 0,$$
$$[X^{i}, g^{jk}] = 0,$$ $$[X^{i}, g_{jk}] = 0.$$
Here it is assumed that $$g^{ij}g_{jk} = \delta^{i}_{k}$$ and that
$$g_{ij}g^{jk} = \delta_{i}^{k}.$$
\vspace{3mm}

The first result that is a direct consequence of these assumptions coupled
with the discrete ordered calculus is the symmetry of the "metric"
coefficients $g^{ij}.$  That is, we shall show that
$$g^{ij}=g^{ji}.$$
\vspace{3mm}

\noindent
{\bf Lemma 3.} $g^{ij}=g^{ji}.$
\vspace{2mm}

\noindent
{\bf Proof.}
$$g^{ij}$$
$$= [X^{i}, \dot{X^{j}}] $$
$$= [X^{i}, [X^{j},J] ]$$
$$= -[ J , [X^{i}, X^{j}] ] - [X^{j} , [J , X^{i}] ]$$
$$= -[J,0] + [X^{j} , [X^{i}, J] ]$$
$$= [X^{j} , [X^{i}, J] ]$$
$$= g^{ji}.$$
\vspace{3mm}

A stream of consequences then follow by differentiating both sides of the
equation $$g^{ij}= [X^{i}, \dot{X^{j}}] $$  where
$$\dot{F} = \dot{X^{j}} \partial_{j} F$$
 and it is understood that $$\partial_{j} F = [F, \dot{X_{j}}] =[F, g_{jk}
\dot{X^{k}}]$$ for any function $F$ of the variables $X^{k}$ and their
derivatives $\dot{X^{k}}.$  In particular, the Levi-Civita connection
$$\Gamma^{ijk} =(1/2)(\partial^{j}g^{ki}
+\partial^{k}g^{ij}-\partial^{i}g^{jk})$$  associated with the $g^{ij}$
comes up almost at once from the differentiation process described above.
One finds that
$$D^{2}X^{i} = G^{i} -
g^{ir}g^{js}F_{rs}\dot{X_{j}}-\Gamma^{ijk}\dot{X_{j}} \dot{X_{k}}$$
where $F_{rs} = [\dot{X_{r}}, \dot{X_{s}}].$ It follows from the Jacobi
identity that $F_{rs}$ satisfies the equation
$$\partial_{i}F_{jk} + \partial_{j}F_{ki} + \partial_{k}F_{ij} = 0,$$
identifying $F_{ij}$ as a noncommutative analog of a gauge field. In a more
technical sense, $G^{i}$ is a noncommutative analog of a scalar field,
satisfying $$<\partial_{i} G_{j}> = <\partial_{j} G_{i}>$$
where the brackets around this equation indicate an analog of the Weyl
ordering for operator products.  The details of this calculation can be
found in \cite{KN:DG}.
\vspace{3mm}

This brief technical description of the equations for a noncommutative
particle in a metric field illustrates well the role of the background of
discrete time in this theory. In terms of the backgound time the metric
coefficients are not constant. It is through this variation that the
spacetime derivatives of the theory are articulated. Thus we are in this
way producing the beginnings of a theory of spacetime based on a background
process. The background is a process with its own form of discrete time,
but no spacetime structure as we know and observe it. Our observation of
spacetime structure appears as a rough (commutative) approximation to the
processes described as consequences of the basic noncommutative equations
of the discrete ordered calculus.
\vspace{3mm}

\section{Poisson Brackets and Commutator Brackets}

Dirac \cite{Dirac} introduced a fundamental relationship between quantum
mechanics and classical mechanics that is summarized by the maxim {\em
replace Poisson brackets by commutator brackets.} Recall that the Poisson
bracket $\{ A, B\}$ is defined by the formula

$$\{ A, B \} = (\partial A/ \partial q) (\partial B/ \partial p) -
(\partial A/ \partial p) (\partial B/ \partial q),$$

\noindent
where $q$ and $p$ denote classical position and momentum variables respectively.
\vspace{3mm}

In our version of discrete physics the noncommuting variables are functions
of discrete time, with a DOC derivative $D$ as described in the previous
section. Since $DX=XJ -JX = [X,J]$ is itself a commutator, it follows that
$$D([A,B]) = [DA,B] + [A, DB]$$ for any expressions $A$, $B$ in our ring
$R$. A corresponding Leibniz rule for Poisson brackets would read
$$(d/dt) \{ A, B \} = \{ dA/dt , B \} + \{ A, dB/dt \}.$$
\vspace{3mm}

\noindent
However, here there is an easily verified exact formula: $$(d/dt) \{ A, B
\} = \{ dA/dt , B \} + \{ A, dB/dt \} - \{ A, B \}(\partial \dot{q} /
\partial q + \partial \dot{p} / \partial p).$$

\noindent
This means that the Leibniz formula will hold for the Poisson bracket
exactly when
$$(\partial \dot{q} / \partial q + \partial \dot{p} / \partial p)=0.$$
\vspace{3mm}

\noindent
This is an integrability condition that will be satisfied if $p$ and $q$
satisfy Hamilton's equations
$$ \dot{q} = \partial H / \partial p ,$$ $$ \dot{p} = - \partial H /
\partial q. $$
\vspace{3mm}

\noindent
This, of course, means that $q$ and $p$ are following a principle of least
action with respect to the Hamiltonian $H$. Thus we can interpret the {\em
fact} $D([A,B]) = [DA,B] + [A, DB]$ in the discrete context as an analog of
the principle of least action. Taking the discrete context as fundamental,
we say that Hamilton's equations are {\em motivated} by the presence of the
Leibniz rule for the discrete derivative of a commutator. The classical
laws are obtained by following Dirac's maxim in the opposite direction!
Classical physics is produced by following the correspondence principle
upwards from the discrete.
\vspace{3mm}

Taking the  last paragraph seriously, we must reevaluate the meaning of
Dirac's maxim.  The meaning of quantization has long been a  basic mystery
of quantum mechanics. By traversing this territory in reverse, starting
from the noncommutative world, we begin these questions anew.
\vspace{3mm}

\section{Scalar Variables, Chaos and Representations of the Discrete
Ordered Calculus}
The purpose of this short section is to point out the inherent
noncommutativity of the operators in any theory based on the discrete
ordered calculus. It is natural to hope for actual scalar variables in the
course of articulating a theory based on DOC.
\vspace{3mm}

Consider the equation $[X,DX]=Jk$ where $k$ is a constant. This reads
$$J(X'X'-2X'X + XX) =Jk$$ and hence we may consider solutions to the equation
$$(X'X'-2X'X + XX)=k.$$ If X and X' commute then this becomes $$(X- X')^{2}=k$$
with the solution $$X' = X \pm k^{1/2}.$$ For some problems it may be
sufficient to
consider the situation where the variables are successively incremented or
decremented
by a constant.
\vspace{3mm}

The problems arise when we go to more than one variable. For example,
consider the equation $$[X_{i},DX_{j}] = J k\delta_{ij}$$where $i$ and $j$
range from
$1$ to $3$. Then for $i \neq j$ we have $$[X_{i},DX_{j}] = 0.$$ Let $X_{i}
= A$ and
$X_{j} = B.$ Then this equation reads $$AJ(B-B') -J(B-B')A = 0.$$ Hence
$$A'(B-B') - (B-B')A =0.$$ Thus if $A$ and $B$ commute, we conclude that
$(A'-A)(B'-B)=0.$ Unfortunately, this contradicts the equations $[A,DA] =
Jk$ and
$[B,DB] = Jk$ that are given by our assumptions, except in the case where $k=0.$
This analysis shows that noncommutativity of the dynamical variables in
theories based
on the discrete ordered calculus is a part of life.
\vspace{3mm}

{\bf Example.} Noncommutativity can have a
scalar source.   For example, suppose that $X=DT$ where $T$ and $T'$ are
commuting scalars.
Consider the equation $$[X,DX]=J^{2}k$$ where $k$ is a commuting scalar
constant.
Then we have $[DT,DDT] = J^{2}k.$ Let $$\Delta = T'-T$$ and note that
$\Delta$ is also a
commuting scalar. Then $DT=J \Delta,$ and therefore $$[DT,DDT] =
J^{2}(\Delta'' (\Delta' - \Delta) - (\Delta'' - \Delta') \Delta).$$
Hence the equation $[X,DX]=J^{2}k$ translates into $$\Delta'' (\Delta' -
\Delta) - (\Delta'' - \Delta') \Delta = k,$$ whence
$$\Delta'' = (k - \Delta \Delta')/(\Delta' - 2 \Delta).$$ This recursion
relation for $\Delta$ and its time series has remarkable properties.
For a fixed non-zero value of $k$, the recursion is highly sensitive to initial
conditions, with regions that give rise to bounded oscillations and other
regions that
give rise to unbounded oscillations. There are boundary values in the initial conditions
where the system undergoes chaotic transition between bounded and unbounded
behaviour.
\vspace{3mm}

We are  investigating this method (of letting $X_{i} = D^{n}T_{i}$ for some
$n$ where $T_{i}$ and $T'_{j}$ are commuting scalars) for producing a
system of noncommuting extrinsic dynamical variables with an underlying
scalar structure. If this idea is correct, then there will emerge a picture
of noncommutative discrete physics based on DOC as a global description
occurring over a substrate of discrete chaotic dynamics.
\vspace{3mm}

There are other possibilites for the direct representation of the discrete
noncommutative dynamics.  There may be matrix representations of these
theories over finite fields, the simplest cases being modular number
systems with prime modulus. This subject will be taken up in a future
publication.
\vspace{3mm}

\section{Discussion on $q$-Deformation}
The direct relation between the content of local physical descriptions
based on the DOC calculus and more global considerations are a matter of
speculation.   One strong hint is contained in the properties of the
discrete derivative that has the form $$D_{q}f(x) = (f(qx) -
f(x))/(qx-x).$$ The classical derivative occurs in the limit as $q$
approaches one.
\vspace{3mm}

In the setting of $q$ not equal to one, the derivative $D_{q}$ is directly
related to fundamental noncommutativity.   Consider variables $x$ and $y$
such that $yx=qxy$ where $q$ is a commuting scalar. Then the expansion of
$(x+y)^n$ generates a $q$-binomial theorem with $q$-choice coefficients
composed in $q$-factorials of $q$-integers $[n]_{q}$ where $$[n]_{q} = 1 +
q + q^2 + ... + q^{(n-1)}.$$ The derivative $D_{q}$ is directly related to
the $q$-integers via the formula $$D_{q}(x^{n}) = [n]_{q} x^{n-1}.$$
\vspace{3mm}

In the context of this paper, we have considered discrete derivatives in
the form $$d_{\Delta}f(x) = (f(x+\Delta) - f(x))/\Delta.$$ This will
convert to the $q$-derivative if $x+ \Delta = qx$.  Thus we need
$$q = (x+ \Delta)/x.$$  This means that a direct translation from DOC to
$q$-derivations could be effected if we allowed $q$ to vary as a function
of $x$ and introduced the temporal operator $J$  into the calculus of
$q$-derivatives.
\vspace{3mm}

In general, many $q$-deformed structures such as the quantum groups
associated with the classical Lie algebras appear to be entwined with the
discretization inherent in $D_{q}.$ The quantum groups have turned out to
be deeply connected with topological amplitudes for networks describing
knots and three dimensional spaces. (See the next section of this paper.)
The analog for the quantum groups in dimension four is being sought. If
there is a connection between the local and the global parts of our essay
it may lie in hidden connections between discretization and quantum groups.
Clearly there is much work to be done in this field.
\vspace{3mm}

There is a clue about the meaning of the operator $J$ ($DF = [F,J]$ in the
discrete ordered calculus)  in the context of quantum groups.  Quantum
groups are Hopf algebras. A quantum group such as $G=U_q(SU(2))$ is
actually an algebra over a field $k$ with an antipode $$S:G \longrightarrow
G$$ and a coproduct $$\Delta: G \longrightarrow G \otimes G$$ , a unit $1$
and a couinit $$\epsilon: G \longrightarrow k$$ The coproduct is a map of
algebras.  The antipode is an antimorphism,  $S(xy) = S(y)S(x),$
and generalizes the inverse in a group in the sense that
$\Sigma S(x_{1})x_{2} = \epsilon(x)1$ and   $\Sigma x_{1}S(x_{2}) =
\epsilon(x)1$ where $\Delta(x) = \Sigma  x_{1} \otimes x_{2}.$
\vspace{3mm}

An element $g$ in a quantum group $G$ is said to be a {\em grouplike
element} if $\Delta(g) = g \otimes g$  and $S(g) = g^{-1}.$
In many quantum groups (such as  $G=U_q(SU(2))$) the square of the antipode
is represented via conjugation by a special grouplike element that we shall
denote by $J$.  Thus  $$S^{2}(x) = J^{-1}xJ$$ for all $x$ in $G.$ This
means that it is possible to define the discrete ordered calculus in the con
text of a quantum group $G$ (as above) by taking $J$ to be the special
grouplike element.  Then we have
$$DX = [X,J] = XJ - JX = J(J^{-1}XJ - X) = J(S^{2}(X) -X).$$
Conjugation by the special grouplike element in the quantum group
constitutes the time evolution operator in this algebra.
\vspace{3mm}

There are a number of curious aspects to this use of the discrete ordered
calculus in a quantum group. First of all, it is the case that in some
quantum groups (for example with undeformed classical Lie algebras) the
square of the antipode is equal to the identity mapping. From the point of
view of DOC, time does not exist in these algebras. But in the
$q$-deformations such as $U_q(SU(2))$, the square of the antipode is quite
non-trivial and can serve well as the tick of the clock.  In this way,
$q$-deformations do provide a context for time. In particular, this
suggests that the $q$-deformations of classical spin networks
\cite{Pen:Spin} should be able to accommodate time. A suggestion directly
related to this remark occurs in \cite{Crane2}, and we shall take this up
at the end of the next section of this paper.
\vspace{3mm}

\section {Networks and Discrete Spacetime}
One can consider replacing continuous space (such as Euclidean space with
the usual topology) by a discrete structure of relationships. The geometry
of the Greeks held a discrete web of relationships in the context of
continuous space. That space was not coordinatized in our way, nor was it
held as an infinite aggregate of points. In general topology there is a
wide choice for possible spatial structures (where we mean by a space a
topology on some set).
\vspace{3mm}

Discretization of space and time implicates the replacement of spacetime by
a network, graph or complex that has nodes for the points and edges to
indicate significant relationships among the points.
\vspace{3mm}

Euler's work in the eighteenth century  brought forth the use of abstract
graphs as holders of spatial structure. After Euler it was possible to find
the classification of the Greek regular solids in the the (wider)
classification of the regular graphs on the surface of the sphere. Metric
can disappear into relationship under the topological constraint of Euler's
formula $V-E+F=2$, where $V$ denotes the number of vertices, $E$ the number
of edges and $F$ the number of faces for the connected graph $G$ on the
sphere.
\vspace{3mm}

A network itself can represent an abstract space.  Embeddings of that
network into a given space (such as graphs on the two dimensional sphere)
correspond to global constraints on the structure of the abstract graph.
\vspace{3mm}

Now a new theme arises, motivated by a conjunction of combinatorics and
physics. Imagine labelling the edges of the network from some set of
"colors". These colors can represent the basic states of a physical system,
or they can be an abstract set of distinct markers for purely mathematical
purposes. Once the network is labelled, each vertex is an entity with a
collection of labels incident to it. Let there be given a function that
associates a number (or algebra element) to each such labelled vertex. Call
this number the {\em vertex weight} at that vertex. Let $C$ denote a
specific coloring of the network $N$ and  consider   the product, over all
the vertices of $N$ of the values of the vertex weights. Finally let $Z(N)$
, the {\em amplitude} of the network,  be defined as the summation of  the
product of the vertex weights over all colorings of the net.  $Z(N)$  is
also called the {\em partition function} of the network.
\vspace{3mm}

Amplitudes of this sort are exactly what one computes in finding the
partition function of a physical system or the quantum mechanical amplitude
for a discrete process. In all these cases the network is interwoven with
the algebraic structure of the vertex weights. It is only recently that
topological properties of networks in three dimensional space have come to
be understood in this way
\cite {Kauff:KP}, \cite{Atiyah},\cite{Witten:QFTJP}. This has led to new
information about thetopology of low dimensional spaces, and new
relationships between physics and topology.
\vspace{3mm}

A classical example of such an amplitude was discovered by Roger Penrose
\cite {Penrose} in elucidating special colorings of 3-regular graphs in the
plane. A 3-regular graph $G$ has three edges incident to each vertex. When
embedded in the plane, these edges acquire a specific cyclic order. Three
colors are used. One associates to each vertex the weight $$\sqrt{-1}\
\epsilon_{abc}$$ where $a$,$b$,$c$ denote the edges meeting the vertex in
this cyclic order, and the epsilon is equal to $1$, $-1$ according as the
edges have distinct labels in the given or reverse cyclic order, or $0$ if
there is a repetition of labels. The resulting amplitude counts the number
of ways to color the network with three colors so that three distinct
colors are incident to each vertex. This result is a perspicuous
generalization of the classical four color problem of coloring maps in the
plane with four colors so that adjacent regions receive different colors.
\vspace{3mm}

The Penrose example generalizes to networks whose amplitudes embody
geometrical properties of Euclidean three dimensional space (angles and
their dependence). Geometry begins to emerge in terms of  the  averages of
properties of an abstract and discrete network of relationships.
Topological properties emerge in the same way. The idea of space may
change to the idea of a network with global states and a functor that
associates this network and its states to the more familiar properties that
a classical observer might see.
\vspace{3mm}

\subsection{Remarks on Quantum Mechanics}

We should remark on the basic formalism for amplitudes in quantum
mechanics. The Dirac notation $\langle A|B\rangle$ \cite{Dirac} denotes the
probability amplitude for a transition from $A$ to $B$. Here $A$ and $B$
could be points in space (for the path of a particle), fields (for quantum
field theory), or geometries on spacetime (for quantum gravity). The
probability amplitude is a complex number. The actual probability of an
event is the absolute square of the amplitude. If a complete set of
intermediate states $C_{1}, C_{2},...C_{n}$ is known, then the amplitude
can be expanded to a summation
$$\langle A|B\rangle = \Sigma_{i=1}^{n}\langle A|C_{i}\rangle\langle
C_{i}|B\rangle.$$ This formula follows the formalism of the usual rules for
probability, and it allows for the constructive and destructive
interference of the amplitudes. It is the simplest case of a quantum
network of the form $$A---*---C--- *---B$$ where the colors at $A$ and $B$
are fixed and we run through all choices of colors for for the middle edge.
The vertex weights at the vertices labelled $*$ are $\langle A|C\rangle$
and $\langle C|B\rangle$ respectively. A measurement at the $C$ edge
reduces the big summation to a single value.
\vspace{3mm}

Consider the generalization of the previous example to the graph

$$A---*---C^{1}---*---C^{2}---*--- ... ---*---C^{m}---B$$

With A and B fixed the amplitude for the net is
$$<A|B> = \Sigma_{1 \leq i_{1} \leq ... \leq i_{m} \leq n}
<A|C^{1}_{i_{1}}><C^{2}_{i_{2}}|C^{3}_{i_{3}}>...<C^{m}_{i_{m}}|B>$$

One can think of this as the sum over all the possible paths from $A$ to
$B.$ In fact in the case of a "particle" travelling between two points in
space, this is exactly what must be done to compute an amplitude -
integrate over all the paths between the two points with appropriate
weightings.  In the discrete case this sort of summation makes perfect
sense. In the case of a continuum there is no known way to make rigorous
mathematical sense out of all cases of such integrals. Nevertheless, the
principles of quantum mechanics must be held foremost for physical purposes
and so such "path integrals" and their generalizations to quantum fields
are in constant use by theoretical physicists \cite{Feynman&Hibbs} who take
the point of view that the proof of a technique is in the consistency of
the results with the experiments.  When the observations themselves are
mathematical (such as finding invariants of knots and links), the issue
acquires a new texture.
\vspace{3mm}

Now consider the summation discussed above in the case where $n=2.$ That
is, we shall assume that each $C^{k} $can take two values, call these
values $L$ and $R.$ Furthermore let us suppose that $<L|R> = <R|L>= \surd
\overline {-1}$ while $<L|L>=<R|R>=1.$ The amplitudes that one computes in
this case correspond to solutions to the Dirac equation \cite{Dirac} in one
space variable and one time variable. This example is related to an
observation of  Richard Feynman
\cite{Feynman&Hibbs}. In \cite{KN:Dirac} we give a very elementary
derivation of this result and we show how  these amplitudes give solutions
to the discretized Dirac equation, so everything is really quite exact and
one can understand just what happens in taking the limit to the continuum.
In this example a state of the network consists in a sequence  of choices
of $L$ or $R$. These can be interpreted as choices to move left or right
along the light-cone in a Minkowski plane. It is in summing over such paths
in spacetime that the solution to the Dirac equation appears.  In this
case, time has been introduced into the net by interpreting the sequence of
nodes in the network as a temporal direction.
\vspace{3mm}

Thus one way to incorporate spacetime is to introduce a temporal direction
into the net. At a vertex, one must specify labels of {\em before} and {\em
after} to each edge of the net that is incident to that vertex. If there is
a sufficiently coherent assignment of such local times, then a global time
direction can emerge for the entire network. Networks endowed with temporal
directions have the structure of morphisms in a category where each
morphism points from past to future. A category of quantum networks emerges
equipped with a functor (via the algebra of the vertex weights) to
morphisms of vector spaces and representations of generalized symmetry
groups. Appropriate traces of these morphisms produce the amplitudes.
\vspace{3mm}

Quantum non-locality is built into the network picture. Any observer taking
a measurement in the net has an effect on the global set of states
available for summation and hence affects the possibilities of observations
at all other nodes in the network. By replacing space with a network we
obtain a precursor to spacetime in which quantum mechanics is built into
the initial structure.
\vspace{3mm}

\noindent
{\bf Remark.} A striking parallel to the views expressed in this section
can be found in \cite{Ett}. Concepts of time and category are discussed by
Louis Crane \cite{Crane1}, \cite{Crane2} in relation to topological quantum
field theory. In the case of Crane's work there is a deeper connection with
the methods of this paper, as I shall explain below.
\vspace{3mm}

\subsection{Temporality and the Crane Model for Quantum Gravity}

Crane uses a partition function defined for a triangulated four-manifold.
Let us denote the partition function by $Z(M^{4}, A, B) = <A|B>_{M}$ where
$M^{4}$ is a four-manifold and $A$ and $B$ are (colored - see the next
sentence) three dimensional submanifolds in the boundary of $M$. The
partition function is constructed by summing over all colorings of the
edges of a dual complex to this triangulation from a finite set of colors
that correspond to certain representations of the the quantum group
$U_{q}(SU(2))$ where $q$ is a root of unity.  The sum is over products of
$15J_{q}$ symbols (natural generalizations of the $6J$ symbols in angular
momentum theory) evaluated with respect to the colorings.  The specific
form of the partition function (here written in the case where $A$ and $B$
are empty)  is

$$Z(M^{4}) = N^{v-e} \Sigma_{\lambda} \Pi_{\sigma} dim_{q}(\lambda(\sigma))
\Pi_{\tau} dim^{-1}_{q}(\lambda(\tau)) \Pi_{\zeta}
15J_{q}(\lambda(\zeta)).$$

Here $\lambda$ denotes the labelling function, assigning colors to the
faces and tetrahedra of $M^{4}$ and $v-e$ is the difference of the number
of vertices and the number of edges in $M^{4}.$  Faces are denoted by
$\sigma$, tetrahedra by $\tau$and 4-simplices by $\zeta.$ We refer the
reader to \cite{CKY} for further details.
\vspace{3mm}

In computing $Z(M^{4}, A, B)=<A|B>_{M}$ one fixes the choice of coloration
on the boundary parts $A$ and $B$.  The analog with quantum gravity is that
a colored three manifold $A$ can be regarded as a three manifold with a
choice of (combinatorial) metric. The coloring is the combinatorial
substitute for the metric. In the three manifold case this is quite
specifically so, since the colors can be regarded as affixed to the edges
of the simplices. The color on a given edge is interpreted as  the
generalized distance between the endpoints of the edge. Thus $<A|B>_{M}$ is
a summation over "all possible metrics" on $M^{4}$ that can extend the
given metrics on $A$ and $B$. $<A|B>_{M}$ is an amplitude for the metric
(coloring) on $A$ to evolve in the spacetime $M^{4}$ to the metric
(coloring) on $B$.
\vspace{3mm}

The partition function  $Z(M^{4}, A, B) = <A|B>_{M}$
is a topological invariant of the four manifold $M^{4}$. In particular, if
$A$ and $B$ are empty (a vacuum-vacuum amplitude), then  the Crane-Yetter
invariant, $Z(M^{4})$,  is a function of the signature and Euler
characteristic of the four-manifold \cite{CKY}. On the mathematical side of
the picture this is already significant since it provides a new way to
express the signature of a four-manifold in terms of local combinatorial
data.
\vspace{3mm}

From the point of view of a theory of quantum gravity, $Z(M^{4}, A, B) =
<A|B>_{M}$, as we have described it so far, is lacking in a notion of time
and dynamical evolution on the four manifold $M^{4}$. One can think of $A$
and $B$ as manifolds at the initial and final times, but we have not yet
described a notion of time within $M^{4}$ itself.
\vspace{3mm}

Crane proposes to introduce time into $M^{4}$ and into the partition
function  $<A|B>_{M}$ by labelling certain three dimensional submanifolds
of $M^{4}$ with special grouplike elements from the quantum group
$U_{q}(SU(2))$ and extending the partition function to include this
labelling.  Movement across such a labelled hypersurface is regarded as one
tick of the clock.  The special grouplike elements act on the
representations in such a way that the partition function can be extended
to include the extra labels.  Then one has the project to understand the
new partition function and its relationship with discrete dynamics for this
model of quantum gravity.
\vspace{3mm}

Lets denote the special grouplike element in the Hopf algebra  $G =
U_{q}(SU(2))$  by the symbol $J.$  Then, as discussed at the end of the
previous section,  one has that the square of the antipode $S:G
\longrightarrow G$ is given by the formula $S^{2}(x) = J^{-1}xJ.$ This is
the tick of the clock.  The DOC derivative in the quantum group is given by
the formula $DX = [X,J] = J(S^{2}(X) - X).$  I propose to generalize the
discrete ordered calculus on the quantum group to a discrete ordered
calculus on the four manifold $M^{4}$ with its hyperthreespaces labelled
with special grouplikes. This generalised calculus will be a useful tool in
elucidating the dynamics of the Crane model. Much more work needs to be
done in this domain.


\vspace{3mm}





\begin{thebibliography}{99}

\bibitem{Atiyah}
M.F. Atiyah [1990],
{\em The Geometry and Physics of Knots,}
Cambridge University Press.

\bibitem{Penrose}
R. Penrose [1971],
Applications of negative dimensional tensors, In {\em Combinatorial
Mathematics and Its Applications}, edited by D. J. A. Welsh, Academic
Press.

\bibitem{CKY}
Crane, Louis , Kauffman, Louis H., Yetter, David N. [1997],
State sum invariants of 4-manifolds,
{\em Journal of Knot Theory and Its Ramifications}, Vol. 6, No. 2. pp. 177-234.

\bibitem{Connes}
Connes,Alain [1990],
{\em Noncommutative Geometry}
Academic Press.

\bibitem{Crane1}
Crane,Louis [1996],
Clock and category: Is quantum gravity algebraic?,
{\em J. Math. Phys.} {\bf 36} (11), November (1996), pp. 6180-6193.

\bibitem{Crane2}
Crane, Louis [1997],
A proposal for the quantum theory of gravity,
(preprint).

\bibitem{Dirac}
Dirac, P.A.M. [1968],
{\em Principles of Quantum Mechanics,}
Oxford University Press.

\bibitem{Ett}
T. Etter,
Process, System, Causality and Quantum Mechanics', {\em
International Journal of General Systems}, issue edited by K. Bowden on
{\bf General systems and the Emergence of Physical Structure from
Information Theory}, (in press).

\bibitem{Fkin}
E. Fredkin [1990],
Digital Mechanics,
{\em Physica D} {\bf 45}, pp. 254-270.

\bibitem{Smolin:QG}
Ashtekar,Abhay, Rovelli, Carlo and Smolin,Lee [1992], "Weaving a Classical
Geometry with Quantum Threads", {\em Phys. Rev. Lett.}, vol. 69, p. 237.


\bibitem{Feynman&Hibbs}
R.P. Feynman and A.R. Hibbs [1965],
{\em Quantum Mechanics and Path Integrals,} McGraw Hill Book Company.


\bibitem{Kauff:KP}
Kauffman,Louis H.[1991,1994],
{\em Knots and Physics,}
World Scientific Pub.


\bibitem{KN:QEM}
Kauffman,Louis H. and Noyes,H. Pierre [1996],
Discrete Physics and the Derivation of Electromagnetism from the formalism
of Quantum Mechanics,
{\em Proc. of the Royal Soc. Lond. A}, {\bf 452}, pp. 81-95.

\bibitem{KN:Dirac}
Kauffman,Louis H. and Noyes,H. Pierre [1996],
Discrete Physics and the Dirac Equation, {\em Physics Letters A}, 218 ,pp.
139-146.

\bibitem{KN:DG}
Kauffman,Louis H. and Noyes,H.Pierre (In preparation)

\bibitem{Twist}
Kauffman, Louis H.[1996],
Quantum electrodynamic birdtracks,
{\em Twistor Newsletter Number 41}

\bibitem{Pen:Spin}
Penrose,Roger [1971],
Angular Momentum - An Approach to Combinatorial Spacetime, In {\em Quantum
Theory and Beyond},Edited by Ted Bastin Cambridge University Press, pp.
151-180.

\bibitem{Tanimura}
Tanimura,Shogo [1992],
Relativistic generalization and extension to the non-Abelian gauge theory
of Feynman's proof of the Maxwell equations, {\em Annals of Physics, vol.
220}, pp. 229-247.

\bibitem{Witten:QFTJP}
Witten,Edward [1989],
Quantum field Theory and the Jones Polynomial, {\em Commun. Math.
Phys.},vol. 121,pp. 351-399.

\end{thebibliography}
\end{document}